# INTENSITY PROCESS AND COMPENSATOR: A NEW FILTRATION EXPANSION APPROACH AND THE JEULIN–YOR THEOREM

BY XIN GUO AND YAN ZENG

*University of California, Berkeley and Bloomberg LP*

Let $(X_t)_{t\geq 0}$ be a continuous-time, time-homogeneous strong Markov process with possible jumps and let $\tau$ be its first hitting time of a Borel subset of the state space. Suppose $X$ is sampled at random times and suppose also that $X$ has not hit the Borel set by time $t$. What is the intensity process of $\tau$ based on this information?

This question from credit risk encompasses basic mathematical problems concerning the existence of an intensity process and filtration expansions, as well as some conceptual issues for credit risk. By revisiting and extending the famous Jeulin–Yor [*Lecture Notes in Math.* **649** (1978) 78–97] result regarding compensators under a general filtration expansion framework, a novel computation methodology for the intensity process of a stopping time is proposed. En route, an analogous characterization result for martingales of Jacod and Skorohod [*Lecture Notes in Math.* **1583** (1994) 21–35] under local jumping filtration is derived.

**1. Introduction.** This paper is motivated by the following problem from credit risk. Let $(X_t)_{t\geq 0}$ be a continuous-time, time-homogeneous strong Markov process and let $\tau$ be its first hitting time of a Borel subset of the state space. Suppose we take samples of $X$ at random times, and suppose also that $X$ has not hit the Borel set by time $t$. Does the intensity process of $\tau$ exist? And if so, how to calculate it? This question encompasses a basic mathematical problem regarding the existence of intensity process in general and some conceptual and computational issues in credit risk study in particular.

*Intensity* $(\lambda_t)_{t\geq 0}$, *compensator* $(A_t)_{t\geq 0}$, *stopping time* $\tau$. The notion of an intensity process $(\lambda_t)_{t\geq 0}$ of a stopping time is of essential interest in credit risk, especially in the information-based approach pioneered by Duffie and Lando [14]. First, if $\tau$ is the default time of a firm and is a stopping time relative to some filtration $\mathbb{G} = (\mathcal{G}_t)_{t \geq 0}$, then under appropriate technical conditions (such as those in Aven [1]), $\lambda_t$ is the instantaneous likelihood of default at time $t$ conditioned on $\mathcal{G}_t$, the information at time $t$. That is,

$$(1) \qquad \lambda_t = \lim_{h \downarrow 0} \frac{1}{h} P(t < \tau \leq t + h | \mathcal{G}_t) \qquad \text{a.s.}$$

Thus, $\lambda_t$ may be viewed as a first-order approximation of the default probability conditioned on the given information $\mathcal{G}_t$. Second, $(\lambda_t)_{t \geq 0}$ is extremely useful for pricing defaultable derivatives. For example, Lando [28] showed that when the default time is formulated as the first jump time of a Cox process (i.e., a doubly stochastic Poisson process) with intensity $(\lambda_t)_{t \geq 0}$, pricing defaultable zero-coupon bond is almost identical to its default-free counterpart except that the discount factor $r_t$ is replaced by $r_t + \lambda_t$. (See also Duffie, Schroder and Skiadas [15] for the reduced-form pricing approach via $(\lambda_t)_{t \geq 0}$, and Jeanblanc and Rutkowski [25] and Bélanger, Shreve and Wong [6] for generalizations.)

Mathematically, the intensity process $(\lambda_t)_{t \geq 0}$ of a stopping time $\tau$ is associated with the compensator $A$ of $\tau$ with respect to an appropriate filtration. Let $(\Omega, \mathcal{F}, P)$ be a probability space, $\tau$ an arbitrary nonnegative random variable, and $\mathbb{F} = (\mathcal{F}_t)_{t \geq 0}$ any filtration. Since $\tau$ is not necessarily an $\mathbb{F}$-stopping time, $(\lambda_t)_{t \geq 0}$ of $\tau$ is associated with the (possibly) expanded filtration $\mathbb{G} = (\mathcal{G}_t)_{t \geq 0}$ of $\mathbb{F}$, where $\tau$ is a $\mathbb{G}$-stopping time. An increasing measurable process $(A_t)_{t \geq 0}$ is called the $\mathbb{G}$-*compensator* of $\tau$ if $A_0 = 0$, $1_{\{\tau \leq t\}} - A_t$ is a $\mathbb{G}$-martingale and $A$ is $\mathbb{G}$-predictable, that is, $A : \mathbb{R}_+ \times \Omega \to [0, \infty]$ is measurable with respect to the $\sigma$-field generated by $\mathbb{G}$-adapted, left-continuous processes. (This is a special case of the well-known Doob–Meyer decomposition where a submartingale is decomposed into a martingale and an increasing predictable process, the latter called the compensator of the submartingale.) The *intensity process* $(\lambda_t)_{t \geq 0}$ of $\tau$ is then defined as the Radon–Nikodym derivative $(dA_t/dt)_{t \geq 0}$, provided that $A$ is a.s. absolutely continuous with respect to the Lebesgue measure. See Brémaud [5], Chapter II, D7, T12 and T13.

### 1.1. Problems and previous work.

*Problem* 1: *Existence and characterization of* $(\lambda_t)_{t \geq 0}$. Given a stopping time $\tau$ with respect to a given filtration, its intensity $(\lambda_t)_{t \geq 0}$ may not exist. A necessary condition for $\tau$ to have an intensity is that $\tau$ be totally inaccessible. For example, the jump times of a Lévy process are totally inaccessible, while the first hitting times of a Brownian motion under the natural filtration are not. However, total inaccessibility is *not* sufficient for the existence of the intensity of $\tau$; see Davis [10] and Giesecke [17], Proposition 6.1, for such an example. As any stopping time can be decomposed into a totally inaccessible part and an accessible part, Zeng [36] recently showed that the totally inaccessible part can be further uniquely decomposed into an absolute continuous part (i.e., with an intensity) and a singular part. Still, the question of *the existence of the intensity for a given stopping time* remains largely open, except in some special cases.

*Problem* 1A: *Computing the intensity* $(\lambda_t)_{t\geq 0}$ *via Meyer's Laplacian approximation.* The existence problem of intensity has been studied in some cases by identifying specific filtrations and stopping times and by explicit calculations. For instance, Dellacherie ([11], Chapter V, T56) calculated the intensity of an arbitrary positive random variable $\tau$ under the natural filtration of $1_{\{\tau \leq t\}}$. Çetin, Jarrow, Protter and Yildirim [7] studied this problem with a filtration generated by the sign of a Brownian motion where the stopping time $\tau$ was chosen to be the first time a Brownian motion doubles in absolute magnitude after remaining below zero for a certain amount of time. This result is based on Azéma's martingale and the excursion theory of Brownian motion. Sezer [33] considered a similar but more general problem with a filtration generated by the process $(R(X_t))_{t\geq 0}$, where $X$ is a one-dimensional diffusion process and $R$ is a Borel function of the form $\sum_{i=1}^{N} i 1_{(x_{i-1}, x_i]}(x)$. By using stochastic integrals with respect to random measures and a martingale representation theorem, she was able to calculate the intensity process of any totally inaccessible stopping time under this particular filtration.

Starting from the more intuitive definition of the intensity process in equation (1), Duffie and Lando [14] considered the first passage time $\tau$ of a one-dimensional diffusion process under a filtration generated by $1_{\{\tau \leq t\}}$ and a discrete, noisy observation of the process at deterministic times. Their approach is based on the Meyer's Laplacian approximation for the compensator of a point process and a dominated convergence type of result by Aven [1]. Following this approach, Guo, Jarrow and Zeng [18, 19] studied and computed $(\lambda_t)_{t\geq 0}$ for first passage times of general Markov processes, including diffusion processes with jumps, and Markov modulated processes. The filtration is of the delayed type, and a minimal filtration expansion was exploited for computation as in [14].

Despite its intuitive appeal, the Meyer's Laplacian approximation approach for computing $(\lambda_t)_{t\geq 0}$ has a serious drawback. First, it is model-specific: the technical conditions for Aven's theorem have to be verified case by case after the first step computation. Second, if one considers the very simple example of the first passage time of a Brownian motion under its natural filtration, then the technical condition of Aven's theorem will be violated, rendering the entire approach invalid. Thus, this approach may give a wrong "intensity" when it actually does not exist. More importantly, Aven's theorem is not suitable for computing $(\lambda_t)_{t\geq 0}$ when it is path dependent. For example, consider a continuous process $X$ continuously observed except on an interval $[S,T]$ which is *random and nonlinearly ordered*, that is, $0 < P(S < T) < 1$. Then $(\lambda_t)_{t\geq 0}$ exists on $(S,T)$, but not outside of $(S,T)$. Moreover, a direct computation using Aven's theorem will usually involve a nontrivial task of finding a dominating process for the approximation.

*Problem* 1B: *Computing the compensator* $(A_t)_{t\geq 0}$ *via the Jeulin–Yor theorem.* Under a different progressive filtration expansion approach, Jeulin and Yor [27] suggested that calculating the compensator (and its intensity) of $1_{\{\tau\leq t\}}$ under the enlarged filtration $\mathbb{G}$ can be translated into computing the compensator of $Z_t = E\{1_{\{\tau>t\}}|\mathcal{F}_t\}$ under the original filtration $\mathbb{F}$, since $Z$ is an $\mathbb{F}$-supermartingale. Unfortunately, the latter problem is not necessarily easier than the former. And [16] considered the very special case when $Z$ is continuous and decreasing, which is also known as satisfying the $G$-hypothesis. In this case, the $\mathbb{F}$-compensator of $Z$ is simply $Z_0 - Z$ and the usual "grad-log" expression $\frac{1}{Z_t}\,dZ_t$ may be analyzed through algebraic manipulation. However, it is not clear how this approach works for general cases when $Z'(t)$ may not even exist and the grad-log expression no longer holds.

*Problem* 2: *Consistency of filtration expansions.* In the course of studying $(\lambda_t)_{t\geq 0}$ via different computation approaches comes the consistency problem of different filtration expansions.

For any nonnegative random variable $\tau$ and a given filtration $\mathbb{F}=(\mathcal{F}_t)_{t\geq 0}$ for which $\tau$ is not a stopping time, $(\lambda_t)_{t\geq 0}$ of $\tau$ is associated with the expanded filtration $\mathbb{G}=(\mathcal{G}_t)_{t\geq 0}$ of $\mathbb{F}$, where $\tau$ is a $\mathbb{G}$-stopping time. As aforementioned, there is more than one way to expand the filtration $\mathbb{F}$: [14, 18, 19] assumed the minimal filtration expansion, and [27] used more than one filtration expansion including the well-known progressive filtration expansion. The progressive filtration expansion is different from the minimal one, and has been recognized as inappropriate from a modeling perspective, as it requires including the original filtration up to $\infty$. (See Section 2.1 for more discussions.) It is thus of both mathematical and applied interests to understand the implication of these different filtration expansions, and in general, their impact on the default intensity and compensator.

1.2. *Our approach and main results.* In this paper Problem 2 is investigated by revisiting and extending the famous and beautiful Jeulin–Yor theorem [27] under a general filtration expansion framework. This generalization resolves the consistency issue of filtration expansions, thus, unifies the approaches and results in [14, 18, 19] and [16].

Based on this generalization and inspired by earlier work of Elliott, Jeanblanc and Yor [16], a different approach to study and compute $(\lambda_t)_{t\geq 0}$ is proposed. And a class of $(\tau, \mathbb{F})$ is identified to obtain an explicit formula for $(\lambda_t)_{t\geq 0}$ of $\tau$ under the expanded filtration $\mathbb{G}$. To overcome the technical difficulty without the $G$-hypothesis, the notion of local jumping filtrations is introduced, and an analogous characterization result for martingales of Jacod and Skorohod [24] under this filtration is established. These ensure the finite variation property of $Z_t = E\{1_{\{\tau > t\}} | \mathcal{F}_t\}$ [hence, the existence of $Z'(t)$] on stochastic intervals. As an illustration of this method, special cases are presented to reproduce the work in [16] when $Z$ is continuous and decreasing, as well as the results in [14, 18, 19, 26].

Finally, for the very special case of the full filtration where filtration expansion is unnecessary, we provide explicit expressions for $(\lambda_t)_{t\geq 0}$ in terms of Lévy systems for Hunt processes.

*Main results.* Throughout the paper, $(\Omega, \mathcal{F}, (\mathcal{F}_t)_{t\geq 0}, P)$ is a filtered probability space where $\mathbb{F} = (\mathcal{F}_t)_{t\geq 0}$ satisfies the usual hypotheses and $\mathcal{F}$ is any $\sigma$-field containing $\mathcal{F}_\infty = \bigvee_{t\geq 0} \mathcal{F}_t$. $\tau$ is a positive random variable, and $Z_t = E\{1_{\{\tau > t\}} | \mathcal{F}_t\}$. Moreover, $\mathbb{G}$ is any filtration expansion of $\mathbb{F}$ such that $\tau$ is a $\mathbb{G}$-stopping time and

$$\mathcal{G}_t \cap \{t < \tau\} = \mathcal{F}_t \cap \{t < \tau\}. \tag{2}$$

THEOREM 1.1 (Jeulin–Yor: Extension). *Given $\mathbb{F}$, $\tau$ and $\mathbb{G}$, the $\mathbb{G}$-compensator of $\tau$ is*

$$\int_0^{t\wedge\tau} \frac{1}{Z_{s-}} dA_s, \tag{3}$$

*where $A$ is the $\mathbb{F}$-compensator of $Z$.*

Note that $\mathbb{G}$ here can be any filtration expansion satisfying equation (2), instead of a progressive one used in the original statement of [27]. In fact, the progressive and the minimal filtration expansions and the filtration expansions in [27], page 78, and [25], Proposition 5.1, all satisfy equation (2), and hence are special cases. Therefore, this extended Jeulin–Yor theorem immediately solves Problem 2. That is, the following:

COROLLARY 1.1. *Compensators of $\tau$ under the progressive and the minimal expansions of $\mathbb{F}$ are identical.*

Though intuitively clear, to the best of our knowledge we are not aware of any prior mathematical result explicitly stating the consistency of compensators/intensities under different filtration expansions.

In addition, Theorem 1.1 enables us to solve Problem 1 for a class of $(\tau, \mathbb{F})$ for which explicit formula of $(\lambda_t)_{t\geq 0}$ is derived.

More precisely, suppose $(\Omega, \mathcal{F}, (\mathcal{F}_t^X)_{t\geq 0}, X, (P_x)_{x\in E}, \theta)$ is a continuous-time, time-homogeneous strong Markov process, such that $\mathbb{F} = (\mathcal{F}_t)_{t\geq 0}$ is a sub-filtration of the completed natural filtration $\mathbb{F}^X = (\mathcal{F}_t^X)_{t\geq 0}$ of $X$, and $\theta_{(\cdot)}$ is a shift operator. Assume $\tau$ is any positive random variable such that the following assumptions hold:

ASSUMPTION (A). (The first time of $X$ hitting a Borel subset of the state space satisfies this condition.) $\tau > t \Leftrightarrow \tau > s$ and $\tau \circ \theta_s > t - s$ ($\forall t, s \geq 0, t > s$).

ASSUMPTION (B). $S$ and $T$ are a pair of $\mathbb{F}$-stopping times so that $X_S \in \mathcal{F}_S$, and $\{S < T\} \subset \{T < \infty\}$. And $\mathbb{F}$ jumps locally from $S$ to $T$, that is,

$$\mathcal{F}_t \cap \{S \leq t < T\} = \mathcal{F}_S \cap \{S \leq t < T\}. \tag{4}$$

ASSUMPTION (C). (See Section 2.2 for the intuition and motivation of this assumption.) $\exists$ some random variables $V_1 \in \mathcal{F}_\infty^X = \sigma(X_t | t \geq 0)$ and $V_2 \in \mathcal{F}_S$, such that $T - S = g(V_1 \circ \theta_S, V_2)$ on $\{S < T\}$, for some Borel function $g$.

Then, we have the following.

THEOREM 1.2. *Given an expanded filtration $\mathbb{G}$ of $\mathbb{F}$ satisfying equation (2) and Assumptions (A), (B) and (C), let $f$ be the following:*

$$f(x, z, t) = \begin{cases} P_x(\tau > t | g(V_1, z) > t), & \text{if } P_x(g(V_1, z) > t) > 0, \\ 0, & \text{if } P_x(g(V_1, z) > t) = 0. \end{cases}$$

*If $f(X_{S(\omega)}(\omega), V_2(\omega), t)$ is absolutely continuous with respect to the Lebesgue measure on $(0, T(\omega) - S(\omega))$ for a.s. $\omega$, then the $\mathbb{G}$-compensator of $\tau$ (i.e., the $\mathbb{G}$-compensator $\widetilde{N}$ of $N_t = 1_{\{\tau \leq t\}}$) satisfies*

$$\begin{aligned}
&1_{(S,T]} d\widetilde{N}_t \\
&= 1_{\{t < \tau\}} 1_{(S,T]} \Bigg( -\frac{f'(X_S, V_2, t - S)}{f(X_S, V_2, t - S)} dt \\
&\quad + \frac{h(X_S, V_2, t - S)}{f(X_S, V_2, t - S)} \frac{P_{X_S}(g(V_1, z) \in dt - a)}{P_{X_S}(g(V_1, z) \geq t - a)} \bigg|_{a=S, z=V_2} \Bigg),
\end{aligned} \tag{5}$$

*where $f'$ is the derivative of $f$ with respect to $t$ and $h(x, z, t) = f(x, z, t-) - P_x(\tau > t | g(V_1, z) = t)$. Here $P_x(\tau > t | V)$ is a Borel function of $V$ with $P_x(\tau > t | V = t)$ taking value at $t$.*

Finally, under the completed natural filtration of a Hunt process (i.e., a càdlàg strong Markov process with totally inaccessible jump times), a stopping time $\tau$ is totally inaccessible if and only if this process $X$ has a jump at $\tau$ a.s. on $\{\tau < \infty\}$ (see [29], pages 111–116). In this case, the explicit expression of $(\lambda_t)_{t\geq 0}$ for the associated stopping time is given by the Lévy system for the Hunt process.

THEOREM 1.3. *Let $X$ be a Hunt process with Lévy system $(U, K)$. For simplicity, we assume $X$ has infinite lifetime. Let $D$ be a Borel subset of the state space, and $\tau = \inf\{t > 0 | X_t \in D\}$. Suppose $X$ has finitely many jumps on every bounded interval. (For a Hunt process $X$, this condition is equivalent to having a strictly positive first jump time.) Then, under the completed natural filtration $(\mathcal{F}_t^X)_{t\geq 0}$ of $X$, the compensator of $\tau_\Lambda$ is*

$$\int_0^{t\wedge\tau_\Lambda} dU_s \int (1_D(y) + 1_{D^c}(y) P_y(\tau = 0)) K(X_s, dy).$$

*[$\tau_\Lambda$ is the totally inaccessible part of $\tau$ such that*

$$\tau_\Lambda(\omega) = \begin{cases} \tau(\omega), & \text{if } \omega \in \Lambda, \\ \infty, & \text{if } \omega \notin \Lambda, \end{cases}$$

*where $\Lambda = \{\tau = \infty\} \cup \{\tau < \infty, X_\tau \neq X_{\tau-}\}$].*

Here $K(x, dy)$ represents intuitively the expected number per unit time of the jumps $X$ makes from $x$ to $dy$, and the additive functional $U$ is the internal clock for $X$. If $X$ is a continuous-time, time-homogeneous Markov chain, $K(x, dy)$ is its generator matrix and $U_t = t$. If $X$ is a Lévy process, $K(x, dy) = \nu(dy - x)$, where $\nu$ is the Lévy measure and $U_t = t$. (See [21], Section 2, [34] and [32], Chapter VI.)

*Organization of the paper.* Section 2 provides detailed proofs of the main results, as well as discussion and comparison with the existing literature. Section 3 illustrates some old and new explicit forms of $(\lambda_t)_{t\geq 0}$ using our approach. The Appendix is a self-contained treatment of an extended version of the Jeulin–Yor formula.

## 2. Proofs and discussions of main results.

2.1. *Discussions on Theorem* 1.1. The original statement of Theorem 1.1, from the insightful work of Jeulin and Yor [27], was built under a progressive filtration expansion framework.

A progressive expansion $\mathbb{G} = (\mathcal{G}_t)_{t\geq 0}$ in the classical filtration expansion theory (see, e.g., [35]) is defined as follows. Given $(\Omega, \mathcal{F}, (\mathcal{F}_t)_{t\geq 0}, P)$ and a positive random variable $\tau$, define

$$\mathcal{G}_t := \{B \in \mathcal{G}_\infty | \exists B_t \in \mathcal{F}_t, B \cap \{t < \tau\} = B_t \cap \{t < \tau\}\},$$

where $\mathcal{G}_\infty = \mathcal{F}_\infty \vee \sigma(\tau)$.

Note that $\mathcal{F}_\infty \cap \{\tau \leq t\} \subset \mathcal{G}_t$, implying intuitively that the expanded $\sigma$-field $\mathcal{G}_t$ for $t \geq \tau$ includes all the information from the unexpanded filtration $\mathbb{F} = (\mathcal{F}_t)_{t \geq 0}$ up to time $\infty$. This seems unrealistic from a modeling perspective.

An alternative is the minimal filtration expansion $\mathbb{G}' = (\mathcal{G}'_t)_{t \geq 0}$ in [11, 14, 18, 19], so that

$$\mathcal{G}'_t = \mathcal{F}_t \vee \sigma(\tau \wedge t) = \mathcal{F}_t \vee \sigma(\{\tau \leq s\} | s \leq t).$$

Clearly, the minimal expansion is the minimal way to expand a filtration to include $\tau$ as a stopping time. And unlike the progressive expansion, it does not require $\mathcal{F}_\infty \cap \{\tau \leq t\} \subset \mathcal{G}'_t$.

Nonetheless, since $\mathcal{G}'_t \subset \mathcal{G}_t$ and $\mathbb{G} = (\mathcal{G}_t)_{t \geq 0}$ is right-continuous, we have

$$\mathcal{F}_t \cap \{t < \tau\} = \mathcal{G}'_t \cap \{t < \tau\} \subset \left( \bigcap_{u > t} \mathcal{G}'_u \cap \{t < \tau\} \right) \subset \left( \bigcap_{u > t} \mathcal{G}_u \cap \{t < \tau\} \right)$$
$$= \mathcal{G}_t \cap \{t < \tau\} = \mathcal{F}_t \cap \{t < \tau\}.$$

Therefore, both $\mathcal{G}_t$ and the right-continuous augmentation of $\mathcal{G}'_t$ coincide with $\mathcal{F}_t$ on $\{t < \tau\}$. It is this critical observation that led us to the formulation of the filtration expansion satisfying equation (2), and to the extended Jeulin–Yor theorem under this general filtration expansion.

Since the original proof by Jeulin and Yor was fairly cryptic, we provide a self-contained proof in the Appendix for completeness.

2.2. *Proof of Theorem* 1.2 *and discussions*.

*Two technical lemmas.* We first need to introduce the notion of local jumping filtration and establish an analogous characterization result for martingales in [24], by adapting their proof in Section 2, pages 22–23.

LEMMA 2.1. *Let $S$ and $T$ be a pair of $\mathbb{F}$-stopping times. If $\mathbb{F}$ jumps locally from $S$ to $T$ according to equation (4), then any uniformly integrable $\mathbb{F}$-martingale $M$ is a.s. of finite variation on $(S, T]$. $[(S, T] := \varnothing$ on $\{S \geq T\}$.] In particular, $Z_t$ is of finite variation on $(S, T]$.*

Another useful lemma can be proved by modifying the technique for jumping filtrations in [20], pages 160–173, to local jumping filtrations.

LEMMA 2.2. *Let $S$ and $T$ be a pair of $\mathbb{F}$-stopping times. If a filtration $\mathbb{F}$ jumps locally from $S$ to $T$ according to equation (4), then for any $\mathbb{F}$-stopping time $\tau$, there exists $\xi \in \mathcal{F}_S$, such that*

$$\tau_{\{\tau < T\}} = \xi_{\{\xi < T\}}.$$

*Proof of Theorem* 1.2. The key is to find the compensator of $Z_t 1_{\{S \leq t < T\}}$ by exploiting the structure of local jumping filtration and the strong Markov property.

First, for the $\mathbb{F}$-optional projection $Z$ of $(1_{\{\tau > t\}})_{t \geq 0}$, by Bayes' formula and Assumption (B),

$$Z_t 1_{\{S \leq t < T\}} = \frac{P(\tau > t, t < T | \mathcal{F}_S)}{P(t < T | \mathcal{F}_S)} 1_{\{S \leq t < T\}}.$$

Here we have taken the convention that when $P(t < T | \mathcal{F}_S) = 0$, $Z_t := 0$, since in this case $P(\tau > t, t < T | \mathcal{F}_S) = 0$. Using Assumptions (A)–(C), the strong Markov property and conditioning on $\mathcal{F}_S^X$, we have on the event $\{S \leq t < T\}$

$$P(t < T | \mathcal{F}_S) = P_{X_S}(g(V_1, z) > r)|_{r=t-S,\, z=V_2}$$

and

$$P(\tau > t, t < T | \mathcal{F}_S) = P(\tau > S, \tau \circ \theta_S > t - S, t - S < g(V_1 \circ \theta_S, V_2) | \mathcal{F}_S)$$
$$= P(\tau > S | \mathcal{F}_S) P_{X_S}(\tau > r, g(V_1, z) > r)|_{r=t-S,\, z=V_2}.$$

So

(6) $\qquad Z_t 1_{\{S \leq t < T\}} = f(X_S, V_2, t - S) P(\tau > S | \mathcal{F}_S) 1_{\{S \leq t < T\}}.$

Next, by Lemma 2.1, $Z$ is of finite variation on $(S, T]$, since $Z$ is an $\mathbb{F}$-supermartingale. If $A$ is the $\mathbb{F}$-compensator of $Z$, then we must have

(7) $\qquad -E\left\{\int_0^\infty H_t 1_{\{S < t \leq T\}} dZ_t\right\} = E\left\{\int_0^\infty H_t 1_{\{S < t \leq T\}} dA_t\right\}$

for any bounded $\mathbb{F}$-predictable process $H$. By the monotone class theorem, in order to find the restriction of $A$ to $(S, T]$, it suffices to consider $H_t = 1_{\{R < t\}}$, where $R$ is an $\mathbb{F}$-stopping time.

By Lemma 2.2, there exists a $\xi \in \mathcal{F}_S$ such that $R_{\{R < T\}} = \xi_{\{\xi < T\}}$. Thus, equation (6) implies

(8)
$$\int_0^\infty H_t 1_{\{S < t \leq T\}} dZ_t$$
$$= \int_0^\infty H_t 1_{\{S < t \leq T\}} f'(X_S, V_2, t - S) P(\tau > S | \mathcal{F}_S) dt + H_T \Delta Z_T 1_{\{S < T\}}.$$

Note that $\{S < T, \xi < T\} \in \mathcal{F}_T$. And we have by first conditioning on $\mathcal{F}_S^X$ and then on $\mathcal{F}_S$, with the strong Markov property and Assumptions (A)–(C), and noting $1_{\{S < T\}} P(t \leq T - S | \mathcal{F}_S^X) = 1_{\{S < T\}} P_{X_S}(t \leq g(V_1, z))|_{z=V_2}$,

$E\{H_T \Delta Z_T 1_{\{S < T\}}\}$
$\qquad = E\{1_{\{S < T, \xi < T\}} [P(\tau > T | \mathcal{F}_T) - f(X_S, V_2, (T - S)-) P(\tau > S | \mathcal{F}_S)]\}$

$$= P(\tau > S, \tau \circ \theta_S > g(V_1 \circ \theta_S, V_2), S < T, \xi - S < g(V_1 \circ \theta_S, V_2))$$

$$- E\{1_{\{S<T, \xi-S<g(V_1\circ\theta_S,V_2)\}}$$

$$\times f(X_0 \circ \theta_S, V_2, g(V_1 \circ \theta_S, V_2)-) P(\tau > S|\mathcal{F}_S)\}$$

$$= E\left\{-1_{\{S<T\}} P(\tau > S|\mathcal{F}_S) \int_{\xi-S}^{\infty} h(X_S, V_2, u) P_{X_S}(g(V_1, z) \in du)|_{z=V_2}\right\}$$

$$= E\left\{-1_{\{S<T\}} P(\tau > S|\mathcal{F}_S) \right.$$

$$\left. \times \int_0^\infty 1_{\{\xi-S<u\leq T-S\}} \frac{h(X_S, V_2, u) P_{X_S}(g(V_1,z) \in du)}{P_{X_S}(g(V_1,z) \geq u)} \bigg|_{z=V_2} \right\}$$

$$= E\left\{-P(\tau > S|\mathcal{F}_S) \right.$$

$$\times \int_0^\infty H_t 1_{\{S<t\leq T\}}$$

$$\left. \times \frac{h(X_S, V_2, t-S) P_{X_S}(g(V_1,z) \in dt-a)}{P_{X_S}(g(V_1,z) \geq t-a)} \bigg|_{a=S, z=V_2} \right\}.$$

One can easily check that $\frac{P_{X_S}(g(V_1,z)\in du)}{P_{X_S}(g(V_1,z)\geq u)}$ in the above equality is well defined by showing that $\{u : Y_u = 0\}$ is a.s. negligible for the regular conditional distribution $P(T \in du|\mathcal{F}_S)$, where $Y_t := P(T \geq t|\mathcal{F}_S)$ and can be chosen as left-continuous.

Finally, combining the above equality with equation (8), we have

$$E\left\{\int_0^\infty H_t 1_{\{S<t\leq T\}} dA_t\right\}$$

(9) $$= -E\left\{\int_0^\infty H_t 1_{\{S<t\leq T\}} dZ_t\right\}$$

$$= E\left\{P(\tau > S|\mathcal{F}_S)\right.$$

$$\times \int_0^\infty H_t 1_{\{S<t\leq T\}}$$

$$\left. \times \frac{h(X_S, V_2, t-S) P_{X_S}(g(V_1,z) \in dt-a)}{P_{X_S}(g(V_1,z) \geq t-a)} \bigg|_{a=S, z=V_2} \right\}$$

$$- E\left\{\int_0^\infty H_t 1_{\{S<t\leq T\}} f'(X_S, V_2, t-S) P(\tau > S|\mathcal{F}_S)\,dt\right\}.$$

Thus, the $\mathbb{F}$-compensator $A$ of $Z$, when restricted to $(S, T]$, satisfies

$$dA_t = P(\tau > S|\mathcal{F}_S)$$

$$\times \left[ -f'(X_S, V_2, t-S) \, dt \right.$$
$$\left. + \frac{h(X_S, V_2, t-S) P_{X_S}(g(V_1, z) \in dt - a)}{P_{X_S}(g(V_1, z) \geq t - a)} \bigg|_{a=S, \, z=V_2} \right].$$

Finally, equation (5) follows from the extended Jeulin–Yor theorem.

*Difference between jumping and local jumping filtrations.* In [24], Jacod and Skorohod studied a filtration $\mathbb{F} = (\mathcal{F}_t)_{t \geq 0}$ generated by a marked point process $(\xi_n, T_n)_{n \geq 1}$, namely, the *jumping filtration*, such that

$$\mathcal{F}_t \cap \{T_n \leq t < T_{n+1}\} = \mathcal{F}_{T_n} \cap \{T_n \leq t < T_{n+1}\}.$$

They proved that any uniformly integrable martingale under this jumping filtration is of finite variation.

It is important to note that a local jumping filtration is not equivalent to a jumping filtration, although the local jumping filtration for the special case of multiple sequences of marked point processes is a jumping filtration. One example relevant to credit risk modeling is given in the following. This example illustrates that the local jumping filtration is critical to resolve both the technical and the computational issues of a path-dependent $(\lambda_t)_{t \geq 0}$ on a random and nonlinearly ordered interval.

EXAMPLE 2.1. *Let $W$ be a one-dimensional Brownian motion starting at level $a$. Let $T_1 = \inf\{t > 0 : W_t = 1\}$ and $T_2 = \inf\{t > 0 : W_t = 2\}$. $a \neq 1, 2$. Clearly, if $a > 2$, then $T_2 < T_1$ a.s.; if $a < 1$, then $T_1 < T_2$ a.s.; if $1 < a < 2$, then $T_1 < T_2$ and $T_1 > T_2$ both have positive probability.*

*Now, consider the process $X_t = W_{t \wedge T_1}$ and $Y_t = \partial 1_{\{T_2 > t\}} + W_t 1_{\{T_2 \leq t\}}$. Here $\partial$ is an "ideal" state different from any point in $\mathbb{R}$. We denote $\mathbb{R} \cup \{\partial\}$ by $E$, the $\sigma$-field $\sigma(X_s, Y_s : s \leq t)$ by $\mathcal{F}_t^0$. We define $\mathcal{F}_t = \bigcap_{u > t} \mathcal{F}_u^0 \vee \mathcal{N}$, where $\mathcal{N}$ is the collection of sets of probability zero. Then $\mathbb{F} = (\mathcal{F}_t)_{t \geq 0}$ satisfies the usual conditions of being complete and right continuous. One can show that this filtration is a local jumping filtration, but not a jumping one.*

*First, $T_1$ and $T_2$ are stopping times under $\mathbb{F}$: $X$ and $Y$ are right continuous with left limit, hence progressively measurable. $T_1$ is the first time $X$ hitting $\{1\}$, and $T_2$ is the first time $Y$ exiting $\{\partial\}$, or equivalently, hitting $\mathbb{R}$.*

*Next, $\mathbb{F}$ is a local jumping filtration, that is,*

$$\mathcal{F}_t \cap \{T_1 \leq t < T_2\} = \mathcal{F}_{T_1} \cap \{T_1 \leq t < T_2\}.$$

PROOF. Clearly, $\mathcal{F}_t \cap \{T_1 \leq t < T_2\} = \mathcal{F}_{t \vee T_1} \cap \{T_1 \leq t < T_2\} \supset \mathcal{F}_{T_1} \cap \{T_1 \leq t < T_2\}$. For the other direction, we first show $\mathcal{F}_t^0 \cap \{T_1 \leq t < T_2\} \subset$ $\mathcal{F}_{T_1} \cap \{T_1 \leq t < T_2\}$: for any $t_1 \leq t_2 \leq \cdots \leq t_n \leq t$ and a Borel subset $B$ of $(E^2)^n$,

$$\{\omega : ((X_{t_1}(\omega), Y_{t_1}(\omega)), \ldots, (X_{t_n}(\omega), Y_{t_n}(\omega))) \in B\} \cap \{T_1 \leq t < T_2\}$$
$$= \{\omega : ((X_{t_1 \wedge T_1}(\omega), \partial), \ldots, (X_{t_n \wedge T_1}(\omega), \partial)) \in B\} \cap \{T_1 \leq t < T_2\}$$
$$\in \mathcal{F}_{T_1} \cap \{T_1 \leq t < T_2\}.$$

By an argument of monotone class theorem, we see $\mathcal{F}_t^0 \cap \{T_1 \leq t < T_2\} \subset \mathcal{F}_{T_1} \cap \{T_1 \leq t < T_2\}$. To extend to the case of $\mathcal{F}_t$, let us suppose $A \in \bigcap_{u > t} \mathcal{F}_u^0$. Then for any $n \geq 1$, there is a set $A_n \in \mathcal{F}_{T_1}$, such that $A \cap \{T_1 \leq t + \frac{1}{n} < T_2\} = A_n \cap \{T_1 \leq t + \frac{1}{n} < T_2\}$. Therefore,

$$A \cap \{T_1 \leq t < T_2\} = A \cap \left( \bigcap_{n=1}^{\infty} \left\{ T_1 \leq t + \frac{1}{n} < T_2 \right\} \right)$$
$$= \bigcap_{n=1}^{\infty} \left( A_n \cap \left\{ T_1 \leq t + \frac{1}{n} < T_2 \right\} \right)$$
$$= \left( \bigcap_{n=1}^{\infty} A_n \right) \cap \{T_1 \leq t < T_2\} \in \mathcal{F}_{T_1} \cap \{T_1 \leq t < T_2\}.$$

Since $A$ is arbitrarily chosen and $\mathcal{F}_{T_1}$ contains all the sets of probability zero, we have proven $\mathcal{F}_t \cap \{T_1 \leq t < T_2\} \subset \mathcal{F}_{T_1} \cap \{T_1 \leq t < T_2\}$. Combined, we have

$$\mathcal{F}_t \cap \{T_1 \leq t < T_2\} = \mathcal{F}_{T_1} \cap \{T_1 \leq t < T_2\}. \qquad \square$$

Intuitively, the filtration is "locally jumping" in the following sense: Suppose $W$ starts from a level between 1 and 2. On the event $\{T_1 < T_2\}$, $W$ will hit level 1 first and then level 2. Accordingly, one observes $W$ from $X$ up to the time $T_1$. Between $(T_1, T_2)$, no observation can be made by the construction of the filtration. After $T_2$, one observes $W$ again from $Y$. So the filtration has a "jump" on $(T_1, T_2)$. Since $(T_1, T_2)$ is a stochastic interval, the "jumping" is "local." Of course, on the event $\{T_2 < T_1\}$, one observes $W$ from $X$ on $(0, T_2)$, from both $X$ and $Y$ on $(T_2, T_1)$, and from $Y$ on $(T_1, \infty)$. So $W$ can be observed all the time and the filtration has no "jump."

*Aven's theorem versus Jeulin–Yor theorem.* The tedious case-by-case technical verification via Aven's theorem is spared in our approach by the use of generalized Jeulin–Yor theorem.

From a computational perspective, our result deals with a general case where $Z$ is only of finite variation, whereas the existing literature assumes $Z_t$ to be monotone and continuous. Most importantly, in our case the intensity may be path-dependent under a local jumping filtration. As such, the extended Jeulin–Yor theorem is a better computation choice than Aven's theorem, which needs to deal with all paths simultaneously and hence unsuitable for path-dependent $(\lambda_t)_{t\geq 0}$.

*Conditions* (A), (B), (C) *in Theorem* 1.2.

(a) Assumption (A) is a generalization of first hitting times of a process. It is particularly useful for the example of [14] in Section 3, where $\tau$ satisfies Assumption (A) but is not the first hitting time of the noisy process.

(b) Assumption (C) is a technical condition to ensure the strong Markov property. Intuitively, it means the time delay $T-S$ is determined in a "nice" way by the future after time $S$ and the history before time $S$. It is satisfied in many important cases, for example, if $t_k$, $t_{k+1}$ are constants with $t_k < t_{k+1}$, $T_n$ and $T_{n+1}$ are, respectively, the $n$th and $(n+1)$st successive jump times or hitting times of $X$, then $S = t_k \vee T_n$ and $T = t_{k+1} \wedge T_{n+1}$ satisfy Assumption (C), as seen from examples in Section 3.

(c) In equation (5), the term $\frac{P_{X_S}(g(V_1,z)\in dt-a)}{P_{X_S}(g(V_1,z)\geq t-a)}|_{a=S,\,z=V_2}$ can be re-written as

$$\frac{P(T\in dt|\mathcal{F}_S)}{P(T\geq t|\mathcal{F}_S)}.$$

This $\mathbb{F}$-intensity of $T$ is a generalization of Jacod [23], Proposition 3.1, under the local jumping filtration.

2.3. *Proof of Theorem* 1.3. First, recall that under the completed natural filtration $(\mathcal{F}_t^X)_{t\geq 0}$, any stopping time $T$ is totally inaccessible if and only if $X$ has a jump at $T$ a.s. on $\{T<\infty\}$ ([29], pages 111–116). This shows the totally inaccessible part of $\tau$ is $\tau_\Lambda$.

For the second part of Theorem 1.3, define $\tau_n^\varepsilon = \inf\{t > \tau_{n-1}^\varepsilon | \rho(X_t, X_{t-}) > \varepsilon\}$ ($\tau_0^\varepsilon = 0$), and $\Lambda^\varepsilon = \bigcup_{n=1}^\infty \{\tau_n^\varepsilon = \tau < \infty\}$. Then $\Lambda = \bigcup_{\varepsilon \in \mathbb{Q}_+} \Lambda^\varepsilon \cup \{\tau = \infty\}$. We note $\tau_n^\varepsilon > 0$ a.s. for $n \geq 1$, so $\{\tau_n^\varepsilon = \tau < \infty\} = \{\tau_n^\varepsilon \leq \tau, \tau_n^\varepsilon < \infty\} \cap A$, where

$A = \{X_{\tau_n^\varepsilon} \in D\} \cup \{X_{\tau_n^\varepsilon} \notin D, (X \circ \theta_{\tau_n^\varepsilon}). \text{ hits D infinitely often near time } 0\}.$

For any $x \in E$, we have $P_x$-a.s.

$$1_{\{\tau_{\Lambda^\varepsilon}\leq t\}} = \sum_{n=1}^\infty 1_{\{\tau_n^\varepsilon = \tau_\Lambda < \infty\}} 1_{\{\tau_n^\varepsilon \leq t\}}$$

$$= \sum_{n=1}^\infty P_x(\tau_n^\varepsilon = \tau_\Lambda < \infty | \mathcal{F}_{\tau_n^\varepsilon}^X) 1_{\{\tau_n^\varepsilon \leq t\}}$$

$$= \sum_{n=1}^{\infty} P_x(\tau_n^\varepsilon \leq \tau_\Lambda, \tau_n^\varepsilon < \infty, A | \mathcal{F}_{\tau_n^\varepsilon}^X) 1_{\{\tau_n^\varepsilon \leq t\}}$$

$$= \sum_{n=1}^{\infty} [1_{\{X_{\tau_n^\varepsilon} \in D\}} + 1_{\{X_{\tau_n^\varepsilon} \notin D\}} P_{X_{\tau_n^\varepsilon}}(\tau = 0)] 1_{\{\tau_n^\varepsilon \leq t \wedge \tau_\Lambda, \tau_n^\varepsilon < \infty\}}.$$

Let $f(x) = 1_D(x) + 1_{D^c}(x) P_x(\tau = 0)$, then $P_x$-a.s.

$$1_{\{\tau_{\Lambda^\varepsilon} \leq t\}} = \sum_{n=1}^{\infty} f(X_{\tau_n^\varepsilon}) 1_{\{\tau_n^\varepsilon \leq t \wedge \tau_\Lambda, \tau_n^\varepsilon < \infty\}} = \sum_{0 < s \leq t} f(X_s) 1_{\{\rho(X_s, X_{s-}) > \varepsilon\}} 1_{[0, \tau_\Lambda]}(s).$$

For any positive predictable process $H$,

$$E_x \left\{ \int_0^\infty H_t \, d1_{\{\tau_{\Lambda^\varepsilon} \leq t\}} \right\} = E_x \left\{ \sum_{n=1}^\infty H_{\tau_n^\varepsilon} f(X_{\tau_n^\varepsilon}) 1_{\{\tau_n^\varepsilon \leq \tau_\Lambda, \tau_n^\varepsilon < \infty\}} \right\}$$

$$= E_x \left\{ \sum_{0 < s < \infty} f(X_s) H_s 1_{[0, \tau_\Lambda]}(s) 1_{\{\rho(X_s, X_{s-}) > \varepsilon\}} \right\}.$$

By the definition of a Lévy system, it is easy to see for any positive function $Y(t, \omega, x)$ that is $\mathcal{P} \otimes \mathcal{B}(E)$-measurable [$\mathcal{P}$ is the predictable $\sigma$-field and $\mathcal{B}(E)$ is the Borel $\sigma$-field on $E$], we have (see [4])

$$E_x \left\{ \sum_{0 < s < \infty} Y(s, X_s) \right\} = E_x \left\{ \int_0^\infty dU_s \int K(X_s, dy) Y(s, y) \right\}.$$

In particular, by setting $Y(s, y) = f(y) H_s 1_{[0, \tau_\Lambda]}(s) 1_{\{\rho(y, X_{s-}) > \varepsilon\}}$, we get

$$E_x \left\{ \int_0^\infty H_t \, d1_{\{\tau_{\Lambda^\varepsilon} \leq t\}} \right\}$$

$$= E_x \left\{ \int_0^\infty H_s 1_{[0, \tau_\Lambda]}(s) \, dU_s \int K(X_s, dy) f(y) 1_{\{\rho(y, X_{s-}) > \varepsilon\}} \right\}.$$

Note

$$E_x \left\{ \int_0^\infty H_t \, d1_{\{\tau_{\Lambda^\varepsilon} \leq t\}} \right\} = E_x \left\{ 1_{\Lambda^\varepsilon} \int_0^\infty H_t \, d1_{\{\tau_\Lambda \leq t\}} \right\},$$

so by letting $\varepsilon \downarrow 0$, we conclude

$$E_x \left\{ \int_0^\infty H_t \, d1_{\{\tau_\Lambda \leq t\}} \right\} = E_x \left\{ \int_0^\infty H_s 1_{[0, \tau_\Lambda]}(s) \, dU_s \int K(X_s, dy) f(y) 1_{\{y \neq X_{s-}\}} \right\}$$

$$= E_x \left\{ \int_0^\infty H_s 1_{[0, \tau_\Lambda]}(s) \, dU_s \int K(X_s, dy) f(y) \right\},$$

since $dU.$ does not charge $\{s > 0 | X_{s-} \neq X_s\}$ and $K(x, \{x\}) = 0$. This shows the compensator of $\tau_\Lambda$ is $\int_0^{t \wedge \tau_\Lambda} dU_s \int f(y) K(X_s, dy)$.

**3. Examples.** In all of the examples below we adopt the minimal filtration expansion if needed, without loss of generality, according to Corollary 1.1.

*Example of Duffie and Lando* [14]. Consider the geometric Brownian motion with noisy observation at a sequence of deterministic times $(t_n)_{n \geq 0}$. More precisely, let $W$ be a standard Brownian motion, $Z_t = Z_0 + mt + \sigma W_t$ and $V_t = e^{Z_t}$. For a constant $V_B > 0$ (this notation is copied from [14] as a default barrier). $\tau$ is defined as $\tau = \inf\{t > 0 | V_t \leq V_B\}$. Let $Y_t = Z_t + U_t$, where $U$ is a Gaussian process independent of $Z$, and for any given $t \in [t_n, t_{n+1})$,

$$\mathcal{H}_t = \sigma(Y_{t_1}, \ldots, Y_{t_n}, 1_{\{\tau \leq s\}} | 0 \leq s \leq t).$$

They calculated the intensity process $\lambda$ of $\tau$ under $\mathbb{H} = (\mathcal{H}_t)_{t \geq 0}$ by

$$\lambda_t = \lim_{h \downarrow 0} \frac{1}{h} P(t < \tau \leq t + h | \mathcal{H}_t).$$

To apply Theorem 1.2, note that if $U$ is another independent Brownian motion, then $Y$ is a strong Markov processes, and $(\mathcal{H}_t)_{t \geq 0}$ is the minimal expansion of a jumping filtration to include $\tau$ as a stopping time. Although $\tau$ is not the first hitting time of $Y$, it still satisfies Assumption (A) as the first hitting time of $V$. For any pair of neighboring observation times $t_k$ and $t_{k+1}$, $t_{k+1} - t_k = (t_{k+1} - t_k) \circ \theta_{t_k}$, $f(Y_{t_k}, t) = P_{Y_{t_k}}(\tau > t)$ $[t \in (0, t_{k+1} - t_k)]$ is differentiable, $P_{Y_S}(t_{k+1} - t_k \in dt - t_k) = \delta_{t_{k+1}}(dt)$ and $h(Y_{t_k}, t_{k+1} - t_k) = P_{Y_{t_k}}(\tau \geq t_{k+1} - t_k) - P_{Y_{t_k}}(\tau > t_{k+1} - t_k) = 0$. Therefore, by Theorem 1.2 and remark (b) on page 12, the intensity process of $\tau$ is

$$1_{\{\tau > t\}} \sum_{k=0}^{\infty} 1_{(t_k, t_{k+1}]} \frac{P_{Y_{t_k}}(\tau \in dt - t_k)}{P_{Y_{t_k}}(\tau > t - t_k)}.$$

Here the true difficulty is finding the distribution function of $\tau$ under the Markovian measures for $Y$. This is a filtering problem and is solved in [14]. In effect, Theorem 1.2 simplifies the calculation by allowing us to use Aven's theorem without checking its technical conditions.

*Example of Guo, Jarrow and Zeng* [18]. Consider a regime switching process $X$ where

$$dX_t = \mu_{\varepsilon(t)} X_t \, dt + \sigma_{\varepsilon(t)} X_t \, dW_t.$$

Here $W$ is a standard Brownian motion and $\varepsilon$ is a continuous-time, time-homogeneous Markov chain independent of $W$, with state space $E = \{0, 1, \ldots, J-1\}$. Let $(t_k)_{k \geq 0}$ be an increasing sequence of deterministic times and $(T_n)_{n \geq 0}$ be the successive jump times of $\varepsilon$. Then $Y = (X, \varepsilon)$ is a strong

Markov process and $T_{n+1} - T_n = T_1 \circ \theta_{T_n}$. The local jumping filtration $\mathbb{F} = (\mathcal{F}_t)_{t \geq 0}$ is generated by the marked point processes $(Y_{t_k}, t_k)_{k \geq 0}$ and $(Y_{T_n}, T_n)_{n \geq 0}$. Let $S = t_k \vee T_n$ and $T = t_{k+1} \wedge T_{n+1}$, then $\mathbb{F}$ jumps locally from $S$ to $T$, and on the event $\{S < T\}$,

$$T - S = T_1 \circ \theta_S \wedge (t_{k+1} - S).$$

So using the notation of Theorem 1.2, $V_1 = T_1$ and $V_2 = t_{k+1} - S$. If $\tau$ is the first passage time of $X$ falling below a level $x$, on the stochastic interval $(0, T - S)$,

$$f(Y_S, V_2, t) = P_{Y_S}(\tau > t | T_1 > t)$$

is differentiable in $t$. On $(0, T - S]$,

$$h(Y_S, V_2, t) = P_{Y_S}(\tau > t | T_1 > t) - P_{Y_S}(\tau > t | T_1 \wedge z = t)|_{z=t_{k+1}-S}.$$

Using Bayes' formula, it is easy to see

$$P_x(\tau > t | T_1 \wedge z) = 1_{\{T_1 \wedge z < z\}} P_x(\tau > t | T_1) + 1_{\{T_1 \wedge z \geq z\}} P_x(\tau > t | T_1 \geq z),$$

$$P_x\text{-a.s.}$$

Therefore,

$$h(Y_S, V_2, t) = [P_{Y_S}(\tau > t | T_1 > t) - 1_{\{t < z\}} P_{Y_S}(\tau > t | T_1 = t)$$
$$- 1_{\{t \geq z\}} P_{Y_S}(\tau > t | T_1 \geq z)]|_{z=t_{k+1}-S}$$
$$= 0.$$

Thus, under the minimal expansion of $\mathbb{F}$, the intensity process of $\tau$ restricted to $(S, T]$ equals ($1_{\{\tau > t\}}$ is omitted for simplicity)

$$-\frac{f'(Y_{t_k \vee T_n}, t - t_k \vee T_n)}{f(Y_{t_k \vee T_n}, t - t_k \vee T_n)} = -\frac{\psi_t(\eta_{\varepsilon(t)}, t - t_k \vee T_n, 1/\sigma_{\varepsilon(t)} \log(x/(X_{t_k \vee T_n})))}{\psi(\eta_{\varepsilon(t)}, t - t_k \vee T_n, 1/\sigma_{\varepsilon(t)} \log(x/(X_{t_k \vee T_n})))},$$

where $\eta_i = \frac{\mu_i}{\sigma_i} - \frac{\sigma_i}{2}$ $(0 \leq i \leq J-1)$, $W_t^{(\eta)} = W_t + \eta t$,

$$\psi(\eta, t, y) = P\left(\inf_{0 \leq s \leq t} W_s^{(\eta)} > y\right) = 1 - \int_0^t \frac{|y|}{\sqrt{2\pi s^3}} e^{-(y-\eta s)^2/(2s)}\, ds$$
(10)
$$\text{for } y < 0,$$

and $\psi_t$ is the derivative of $\psi$ with respect to $t$. This is the same as the result in [18].

*Example of Guo, Jarrow and Zeng* [19]. Consider a jump diffusion process $X$ with

$$X_t = X_0 e^{(\mu - 1/2\sigma^2)t + \sigma W_t} \prod_{0 < s \leq t, \Delta\varepsilon(s) \neq 0} \xi_{\varepsilon(s)}.$$

Here $W$ and $\varepsilon$ are the same as in the previous example. Assume further the generator matrix of $\varepsilon$ is $Q = (q_{ij})_{J \times J}$ with $q_j := \sum_{i \neq j} q_{ji}$. And $(\xi_n)_{n=0}^{J-1}$ is a sequence of random variables, such that $W$, $\varepsilon$ and $(\xi_n)_{n=0}^{J-1}$ are all independent.

Following the notation of the previous example, the marked point processes $(Y_{t_k}, t_k)_{k \geq 0}$ and $(Y_{T_n}, T_n)_{n \geq 0}$ generate a local jumping filtration $\mathbb{F} = (\mathcal{F}_t)_{t \geq 0}$. On $(0, \bar{T} - S)$,

$$f(Y_S, V_2, t) = P_{Y_S}(\tau > t | T_1 > t)$$

is differentiable in $t$. On $(0, T - S]$,

$$\begin{aligned}h(Y_S, V_2, t) &= P_{Y_S}(\tau > t | T_1 > t) - P_{Y_S}(\tau > t | T_1 \wedge z = t)|_{z = t_{k+1} - S} \\ &= 1_{\{t < t_{k+1} - S\}}[P_{Y_S}(\tau > t | T_1 > t) - P_{Y_S}(\tau > t | T_1 = t)] \\ &\quad + 1_{\{t = T - S = t_{k+1} - S\}}[P_{Y_S}(\tau > t | T_1 > t) - P_{Y_S}(\tau > t | T_1 \geq t)] \\ &= 1_{\{T_{n+1} < t_{k+1}\}} P_{Y_S}(\tau = T_1 | T_1 = t).\end{aligned}$$

Moreover, we have

$$\frac{P_{Y_S}(T_1 \wedge z \in dt)}{P_{Y_S}(T_1 \wedge z \geq t)} = 1_{\{t < z\}} \frac{P_{Y_S}(T_1 \in dt)}{P_{Y_S}(T_1 \geq t)} + \delta_z(dt),$$

where $\delta_z(dt)$ is the Dirac measure at $z$. Therefore, under the minimal expansion of $\mathbb{F}$, the intensity process of $\tau$ restricted to $(S, T]$ equals ($1_{\{\tau > t\}}$ is again omitted for simplicity)

$$\begin{aligned}&-\frac{\psi_t(\eta, t - t_k \vee T_n, \frac{1}{\sigma} \log(x/(X_{t_k \vee T_n})))}{\psi(\eta, t - t_k \vee T_n, \frac{1}{\sigma} \log(x/X_{t_k \vee T_n}))} \\ &\quad + 1_{\{T_{n+1} < t_{k+1}\}} \sum_{j \neq \varepsilon(t)} q_{\varepsilon(t)j} \bigg\{ \int_0^1 F_j(dz) \phi\bigg(\eta, t - t_k \vee T_n, \frac{1}{\sigma} \log \frac{x}{X_{t_k \vee T_n}}, \\ &\hspace{9cm} \frac{1}{\sigma} \log \frac{x}{z X_{t_k \vee T_n}}\bigg)\bigg\} \\ &\hspace{5cm} \times \bigg\{\psi\bigg(\eta, t - t_k \vee T_n, \frac{1}{\sigma} \log \frac{x}{X_{t_k \vee T_n}}\bigg)\bigg\}^{-1},\end{aligned}$$

where $\psi$ is defined in equation (10), $\psi_t$ is the derivative of $\psi$ with respect to $t$, and $F_j$ is the distribution function of $\xi_j$. And for $y_1 \leq y_2$ and $W_t^{(\eta)} = W_t + \eta_t$,

$$\phi(\eta, t, y_1, y_2) = P\bigg(\inf_{s \leq t} W_s^{(\eta)} > y_1, W_t^{(\eta)} \leq y_2\bigg)$$

$$= \Phi\left(\frac{y_2 - \eta t}{\sqrt{t}}\right) - \Phi\left(\frac{y_1 - \eta t}{\sqrt{t}}\right)$$
$$- e^{2\eta y_1}\left[\Phi\left(\frac{y_2 - 2y_1 - \eta t}{\sqrt{t}}\right) - \Phi\left(\frac{-y_1 - \eta t}{\sqrt{t}}\right)\right].$$

Here $\eta = \frac{\mu}{\sigma} - \frac{\sigma}{2}$, and $\Phi(x)$ is the distribution of a standard normal random variable. This reproduces the result of [18].

*Example of Elliott, Jeanblanc and Yor* [16]. Given $(\Omega, \mathcal{G}, P)$, a positive random variable $\tau$ and a Brownian motion $B$ on $(\Omega, \mathcal{G}, P)$, consider $\mathbb{F} = (\mathcal{F}_t)_{t \geq 0}$, the completed natural filtration of $B$ for which $\tau$ is not an $\mathbb{F}$-stopping time. Let $\mathbb{G} = (\mathcal{G}_t)_{t \geq 0}$ be the minimal expansion of $\mathbb{F}$ to include $\tau$ as a stopping time. Suppose $Z_t = P(\tau > t | \mathcal{F}_t)$ is continuous and nonincreasing, then the compensator $A$ of $Z$ equals $Z_0 - Z$. By the Jeulin–Yor theorem, the intensity process of $\tau$ equals

$$\lambda_t = -1_{\{\tau > t\}} \frac{Z'_t}{Z_t}.$$

*Example of Guo, Jarrow and Zeng* [19]. Consider the regime switching process as in the previous example, except that $\tau$ is defined as the first time of $Y = (X, \varepsilon)$ hitting the region $D = (-\infty, x) \times \{0\}$ ($x$ is a constant), that is, $\tau := \inf\{t > 0 | X_t < x, \varepsilon(t) = 0\}$. (This choice of $\tau$ was introduced to study the recovery rate process in [19], where notions of bankruptcy and default are differentiated.) Consider $\mathbb{F} = (\mathcal{F}_t)_{t \geq 0}$, the completed natural filtration of $Y$.

Since $X$ is continuous and $\varepsilon$ has a Lévy system $(t, Q)$, $Y = (X, \varepsilon)$ has a Lévy system

$$U_t = t, K((x, i), dy\, dj) = q_{ij}\delta_x(dy) 1_{E - \{i\}}(dj).$$

The totally inaccessible part $\tau_\Lambda$ of $\tau$ has the compensator

$$A_t = \int_0^{t \wedge \tau_\Lambda} ds \int_{-\infty}^{\infty} \sum_{j=0}^{J-1} (1_D(y, j) + 1_{D^c}(y, j) P_{(y,j)}(\tau = 0))$$
$$\times q_{\varepsilon(s)j} \delta_{X_s}(dy) 1_{E - \{\varepsilon(s)\}}(dj).$$

Because for a finite Markov chain each state has nonzero holding time, $P_{(y,j)}(\tau = 0) = 0$ if $j \neq 0$. By the continuity of $X$, $P_{(y,0)}(\tau = 0) = 0$ if $y > x$. Thus,

$$A_t = \int_0^{t \wedge \tau_\Lambda} ds [1_{\{X_s < x\}} + 1_{\{X_s = x\}} P_{(X_s, 0)}(\tau = 0)] 1_{\{\varepsilon(s) \neq 0\}} q_{\varepsilon(s)0}.$$

Hence, $\tau$ has an intensity process

$$1_{\{\varepsilon(t) \neq 0\}} q_{\varepsilon(t)0} 1_{\{X_t \leq x\}}.$$

*Motivating problem* (*Jeanblanc and Valchev* [26]). If the observation times are deterministic, the unexpanded filtration $\mathbb{F} = (\mathcal{F}_t)_{t \geq 0}$ is generated by the marked point process $(X_{t_k}, t_k)_{k \geq 1}$. Similar to the example of [14], by Theorem 1.2, the $\mathbb{G}$-compensator $\widetilde{N}$ of the point process $N_t = 1_{\{\tau \leq t\}}$ has the intensity process

$$1_{\{t < \tau\}} \sum_{k=0}^{\infty} 1_{\{t_k < t \leq t_{k+1}\}} (-1) \frac{f'(X_{t_k}, t - t_k)}{f(X_{t_k}, t - t_k)} \, dt,$$

where $f(x, t) = P_x(\tau > t)$ and $f'$ is the derivative with respect to $t$.

**4. Conclusions.** This paper studied the existence problem of intensity processes and related filtration expansion issues, and provided an alternative computational methodology.

In [24] both necessary and sufficient conditions to characterize martingales under a jumping filtration were given. Our paper derived and exploited the analogous necessity result for the local jumping filtration. However, it remains an interesting mathematical problem to extend their sufficiency result to the local jumping filtration.

Finally, in [2, 3, 9, 22], the notion of "information drift" was introduced and analyzed by tools of Malliavin calculus for the additional information induced by the filtration expansion. Furthermore, this "information drift" was used to identify the additional utility by entropy-related quantities from information theory. Their key technical assumption is that a semi-martingale remains a semi-martingale under filtration expansion. However, in our paper the $\mathbb{G}$-semi-martingale $1_{\{\tau \leq t\}}$ is NOT adapted to $\mathbb{F}$, hence, no longer a semi-martingale under $\mathbb{F}$. (See Stricker's theorem on page 53 of [30] for necessary and sufficient conditions for a semi-martingale to remain a semi-martingale under filtration shrinkage.) This leads to another interesting question: does expanding the original filtration in a minimal way by the nonadapted point process $1_{\{\tau \leq t\}}$ lead any added value/utility in some way? This is an interesting project that goes beyond the scope of this paper (see [8]).

## APPENDIX A: PROOF OF THEOREM 1.1: THE EXTENDED JEULIN–YOR THEOREM

Section A.1 provides the necessary terminology and background. For a more comprehensive presentation of related materials, please refer to [12], Chapter VI or [20], Chapter V.

**A.1. Preliminaries.** Suppose $(\Omega, \mathcal{F}, (\mathcal{F}_t)_{t \geq 0}, P)$ is a complete filtered probability space such that the filtration $\mathbb{F} = (\mathcal{F}_t)_{t \geq 0}$ satisfies the usual hypotheses. The *predictable $\sigma$-field* $\mathcal{P}$ is the $\sigma$-field on $\mathbb{R}_+ \times \Omega$ generated by all the left-continuous processes adapted to $\mathbb{F} = (\mathcal{F}_t)_{t \geq 0}$. Similarly, the *optional $\sigma$-field* $\mathcal{O}$ is the $\sigma$-field generated by all the right-continuous processes adapted to $\mathbb{F} = (\mathcal{F}_t)_{t \geq 0}$. A stopping time $T$ is *predictable* if the process $1_{\{T \leq t\}}$ is $\mathcal{P}$-measurable, or equivalently, there exist stopping times $T_n \uparrow T$ with $T_n < T$ on the set $\{T > 0\}$.

It is important to find "nice versions" of a measurable process, as made precise in the following.

PROPOSITION A.1. *Let $X$ be a positive or bounded $\mathcal{B}(\mathbb{R}_+) \otimes \mathcal{F}$ measurable process. There exist an optional process, denoted by $^oX$, and a predictable process, denoted by $^pX$, such that*

$$E\{X_T 1_{\{T<\infty\}} | \mathcal{F}_T\} = {}^o X_T 1_{\{T<\infty\}} \qquad a.s.,$$
$$E\{X_S 1_{\{S<\infty\}} | \mathcal{F}_S\} = {}^p X_S 1_{\{S<\infty\}} \qquad a.s.$$

*for every stopping time $T$ and every predictable time $S$. $^oX$ and $^pX$ are unique up to indistinguishability and they are called the optional projection and the predictable projection of $X$, respectively.*

REMARK A.2. If $X$ is continuous, $^oX$ and $^pX$ are not necessarily so ([12], page 113, 50(d)).

Another important type of projection processes is the dual projections of an increasing process. (See Dellacherie [11] and Doléans [13].) The main perspective is that an increasing process induces a measure on the product $\sigma$-field $\mathcal{B}(\mathbb{R}_+) \otimes \mathcal{F}$, and the measurability properties of the increasing process are characterized by those of the induced measure.

A $\mathcal{B}(\mathbb{R}_+) \otimes \mathcal{F}$-measurable process $A = (A_t)_{t \geq 0}$ is called a *(raw) increasing process*, if $A$ is right-continuous and increasing. Note $A$ is not necessarily adapted to $\mathbb{F} = (\mathcal{F}_t)_{t \geq 0}$. Define a set function $\mu_A$ on $\mathcal{B}(\mathbb{R}_+) \otimes \mathcal{F}$ by $\mu_A(X) = E\{\int_0^\infty X_t \, dA_t\}$, for any positive measurable process $X \in \mathcal{B}(\mathbb{R}_+) \otimes \mathcal{F}$. Then $\mu_A$ is a measure on $\mathcal{B}(\mathbb{R}_+) \otimes \mathcal{F}$ and since for $T_n = \inf\{t \geq 0 | A_t \geq n\}$, $\bigcup_n [0, T_n) = \mathbb{R}_+ \times \Omega$, $\mu_A$ is $\sigma$-finite.

The following result establishes a one-to-one correspondence between a $\sigma$-finite measure on $\mathcal{B}(\mathbb{R}_+) \otimes \mathcal{F}$ and a raw increasing process.

PROPOSITION A.3. *Suppose $\mu$ is a $\sigma$-finite measure on $\mathcal{B}(\mathbb{R}_+) \otimes \mathcal{F}$. Then there exists a raw increasing process $A$ such that $\mu = \mu_A$ if and only if $\mu(X) = 0$ whenever $X \in \mathcal{B}(\mathbb{R}_+) \otimes \mathcal{F}$ is indistinguishable from 0. Such an $A$ is unique up to indistinguishability.*

Furthermore, the measurability properties of $A$ can be characterized by those of $\mu_A$.

PROPOSITION A.4. *A raw increasing process $A$ is optional (resp. predictable) if and only if $\mu(X) = \mu(^oX)$ [resp. $\mu(X) = \mu(^pX)$] for any positive measurable process $X$.*

REMARK A.5. Since $A$ is right-continuous, $A$ is optional if and only if $A$ is adapted to $\mathbb{F} = (\mathcal{F}_t)_{t\geq 0}$.

Propositions A.3 and A.4 naturally lead to a method of manufacturing "nice" increasing processes from a raw one. To begin with, define via $\mu_A$ a new $\sigma$-finite measure $\mu_A^o$ by letting $\mu_A^o(X) := \mu_A(^oX)$ for any positive measurable process $X$. Since $^o(^oX) = {}^oX$, by Proposition A.4, there exists a unique optional increasing process, denoted by $A^o$, such that $\mu_A^o = \mu_{A^o}$. $A^o$ is called the *dual optional projection of $A$*. Similarly, the *dual predictable projection of $A$* can be defined, denoted by $A^p$, and also is called the *compensator* of $A$.

**A.2. Proof of the extended Jeulin–Yor theorem.** First, recall some general propositions connecting $A$ and its compensator, and a few lemmas.

PROPOSITION A.6. *Let $A$ be a raw increasing process with $E\{A_\infty\} < \infty$, and $B$ a predictable increasing process. Then $B$ is the dual predictable projection of $A$ if and only if $\mu_A$ and $\mu_B$ coincide on the predictable $\sigma$-field, or equivalently, $B_0 = E\{A_0|\mathcal{F}_0\}$ and $^oA - B$ is a uniformly integrable martingale. When $A$ is also adapted, $^oA = A$ and $B$ is the predictable part in the Doob–Meyer decomposition of the submartingale $A$.*

PROPOSITION A.7. *Let $A$ be an adapted, bounded increasing process, then*
$$^p(\Delta A) = \Delta A^p,$$
*where $\Delta A_t := A_t - A_{t-}$ and $\Delta A_t^p := A_t^p - A_{t-}^p$.*

Moreover, according to [27], Lemma 1 or [30], Lemma, page 370, we have the following:

LEMMA A.8. *If $H$ is a $\mathbb{G}$-predictable process, then there exists an $\mathbb{F}$-predictable process $J$ such that $H = J$ on $[0, \tau]$.*

Now, let $Z$ be the $\mathbb{F}$-optional projection of $1_{\{\tau > t\}}$ so that $Z_t = E\{1_{\{\tau > t\}}|\mathcal{F}_t\}$ a.s., then $Z$ is a bounded $\mathbb{F}$-supermartingale, and admits a càdlàg version by [20], Theorems 2.46, 2.47, Chapter II. Moreover, by the Doob–Meyer decomposition theorem, there exists a uniformly integrable $\mathbb{F}$-martingale $M$ and a unique $\mathbb{F}$-predictable increasing process $A$ with $A_0 = 0$ and $E\{A_\infty\} < \infty$, such that $Z = M - A$. Alternatively, $A$ can also be characterized according to Proposition A.6 by the follows.

LEMMA A.9. *A is the $\mathbb{F}$-dual predictable projection of the increasing process $N_t = 1_{\{\tau \leq t\}}$.*

In addition, we have the following:

LEMMA A.10. *The set $\{t | 0 \leq t < \infty, Z_{t-} = 0\}$ is negligible for the measure $dA$.*

PROOF. The key is to show $P(Z_{\tau-} > 0) = 1$. Indeed, assume $P(Z_{\tau-} > 0) = 1$. Since $\{Z_{\cdot -} = 0\}$ is an $\mathbb{F}$-predictable set, we have by Lemma A.9

$$E\left\{\int_0^\infty 1_{\{Z_{t-}=0\}} dA_t\right\} = E\{1_{\{Z_{\tau-}=0\}} 1_{\{\tau < \infty\}}\} = 0.$$

So $\{t | 0 \leq t < \infty, Z_{t-} = 0\}$ is negligible for the measure $dA$.

It remains to show $P(Z_{\tau-} > 0) = 1$. Here we elaborate the argument of [35], Lemma 0. Define $T = \inf\{t \geq 0 | Z_t = 0 \text{ or } Z_{t-} = 0\}$. Then $Z$ vanishes on $[T, \infty)$ according to [31], Proposition 3.4, Chapter II. Since $\tau \leq T$ a.s. from the proof in [30], Theorem 13, the process $(1^Z)_t = \frac{1}{Z_t} 1_{\{t < \tau\}}$ is well defined, and is a $\mathbb{G}$-supermartingale. Indeed,

$$E\{(1^Z)_t\} = E\left\{\frac{1_{\{Z_t \neq 0\}}}{Z_t} 1_{\{t < \tau\}}\right\} = P(Z_t \neq 0) < \infty,$$

and by Bayes' formula,

$$E\{(1^Z)_t | \mathcal{G}_s\} = 1_{\{\tau > s\}} \frac{E\{(1^Z)_t | \mathcal{F}_s\}}{E\{1_{\{\tau > s\}} | \mathcal{F}_s\}} = (1^Z)_s P(Z_t \neq 0 | \mathcal{F}_s) \leq (1^Z)_s.$$

Since $Z$ is càdlàg, $(1^Z)(\omega)$ has a càdlàg sample path for a.s. $\omega$. In particular, $(1^Z)_{\tau-}$ is finite on $\{\tau < \infty\}$. Since $(1^Z)$ is a positive supermartingale, $(1^Z)_\infty := \lim_{t \to \infty} (1^Z)_t$ exists and is integrable. Hence, $(1^Z)_{\tau-} < \infty$ a.s. on $\{\tau = \infty\}$. Combined, we conclude $(1^Z)_{\tau-} < \infty$ a.s. and $P(Z_{\tau-} > 0) = 1$. □

PROOF OF THEOREM 1.1. By Proposition A.6, it suffices to show that for any bounded $\mathbb{G}$-predictable process $H$,

$$E\left\{\int_0^\infty H_t \, dN_t\right\} = E\left\{\int_0^\infty H_t 1_{\{t \leq \tau\}} \frac{dA_t}{Z_{t-}}\right\}.$$

Indeed, for a given $H$, there exists an $\mathbb{F}$-predictable process $J$, such that $H = J$ on $[0, \tau]$. Therefore,

$$E\left\{\int_0^\infty H_t \, dN_t\right\} = E\{H_\tau 1_{\{\tau < \infty\}}\} = E\{J_\tau 1_{\{\tau < \infty\}}\} = E\left\{\int_0^\infty J_t \, dA_t\right\},$$

where the last equality is due to Lemma A.9.

Since $M$ is a uniformly integrable martingale, ${}^pM = M_-$ by the optional sampling theorem. By Proposition A.7, $\Delta A_t = {}^p(\Delta N)_t$. Therefore,

$$Z_{t-} = M_{t-} - A_{t-} = {}^p(M-A)_t + \Delta A_t = {}^p({}^o(1-N)_t) + {}^p(\Delta N)_t = {}^p(1-N_-)_t.$$

Hence,

$$\begin{aligned}
E\bigg\{\int_0^\infty H_t\,dN_t\bigg\} &= E\bigg\{\int_0^\infty \frac{J_t}{Z_{t-}} Z_{t-}\,dA_t\bigg\} = E\bigg\{\int_0^\infty \frac{J_t}{Z_{t-}}{}^p(1-N_-)_t\,dA_t\bigg\} \\
&= E\bigg\{\int_0^\infty \frac{J_t}{Z_{t-}}(1-N_-)_t\,dA_t\bigg\} = E\bigg\{\int_0^\infty \frac{J_t}{Z_{t-}}1_{\{\tau\geq t\}}\,dA_t\bigg\} \\
&= E\bigg\{\int_0^\infty \frac{H_t}{Z_{t-}}1_{\{\tau\geq t\}}\,dA_t\bigg\}. \qquad\square
\end{aligned}$$

**Acknowledgments.** The authors are very grateful to M. Jeanblanc for her inspirational discussions and remarks. We also thank an Associate Editor and the anonymous referees for their constructive suggestions and enlightening remarks.

DEPARTMENT OF INDUSTRIAL ENGINEERING
  AND OPERATIONS RESEARCH
UNIVERSITY OF CALIFORNIA, BERKELEY
4173 ETCHEVERRY HALL
BERKELEY, CALIFORNIA 94720-1777
USA
E-MAIL: xinguo@ieor.berkeley.edu

BLOOMBERG LP
731 LEXINGTON AVENUE
NEW YORK, NEW YORK 10022
USA
E-MAIL: yz44@cornell.edu